\documentclass[a4paper,american,reqno]{amsart}

\newif\ifPreprint \Preprinttrue
\newif\ifSubmission \Submissionfalse

\usepackage{babel}
\usepackage[utf8]{inputenc}
\usepackage[T1]{fontenc}
\usepackage{xspace}
\usepackage{enumitem}
\usepackage{todonotes}
\usepackage{csquotes}
\usepackage{microtype}
\usepackage[style = numeric-comp,
            firstinits = true,
            maxbibnames = 100,
            maxcitenames = 2,
            isbn = false,
            backend = biber]{biblatex}
\usepackage[colorlinks,
            citecolor=blue,
            urlcolor=blue,
            linkcolor=blue]{hyperref}
\usepackage{cleveref}

\makeatletter
\patchcmd{\@settitle}{\uppercasenonmath\@title}{\scshape\large}{}{}
\patchcmd{\@setauthors}{\MakeUppercase}{\scshape\normalsize}{}{}
\makeatother

\vfuzz \hfuzz
\theoremstyle{plain}
\newtheorem{lemma}{Lemma}
\newtheorem{corollary}{Corollary}
\newtheorem{theorem}{Theorem}
\theoremstyle{definition}

\newtheorem{example}{Example}

\newtheorem{remark}{Remark}
\newtheorem{obs}{Observation}

\theoremstyle{remark}

\newcommand{\set}[1]{\{#1\}}
\newcommand{\Set}[1]{\left\{#1\right\}}
\newcommand{\defset}[3][\defsep]{\set{#2#1#3}}
\newcommand{\Defset}[3][\defsep]{\Set{#2#1#3}}
\newcommand{\Unc}{\mathcal{U}}

\newcommand{\st}{\text{s.t.}}

\newcommand{\R}{\mathbb{R}}
\newcommand{\I}{\mathbb{I}}

\newcommand{\tM}{\tilde{M}}

\newcommand{\fa}{\text{ for all }}

\newcommand{\bdt}{\boldsymbol{\cdot}}

\newcommand{\abbr}[1][abbrev]{#1.\xspace}

\newcommand{\eg}{\abbr[e.g]}

\newcommand{\ie}{\abbr[i.e]}

\DeclareMathOperator*{\rk}{rank}

\newcommand{\define}{\mathrel{{\mathop:}{=}}}

\newcommand{\rev}[1]{#1}

\bibliography{adjustable-robust-lcps}

\begin{document}

\title[Affinely Adjustable Robust Linear Complementarity Problems]%
{Affinely Adjustable Robust\\Linear Complementarity Problems}
\author[C. Biefel, F. Liers, J. Rolfes, M. Schmidt]%
{Christian Biefel, Frauke Liers, Jan Rolfes, Martin Schmidt}

\address[C. Biefel, F. Liers, J. Rolfes]{%
  Friedrich-Alexander-Universität Erlangen-Nürnberg,
  Discrete Optimization,
  Cauerstr. 11,
  91058 Erlangen,
  Germany;
  Energie Campus Nürnberg,
  Fürther Str. 250,
  90429 Nürnberg,
  Germany}
\email{\{christian.biefel,frauke.liers,jan.rolfes\}@fau.de}

\address[M. Schmidt]{%
  Trier University,
  Department of Mathematics,
  Universitätsring 15,
  54296 Trier,
  Germany}
\email{martin.schmidt@uni-trier.de}

\date{\today}

\begin{abstract}
  Linear complementarity problems are a powerful tool for
modeling many practically relevant situations such as market
equilibria. They also connect many sub-areas of mathematics like
game theory, optimization, and matrix theory.
Despite their close relation to optimization, the
protection of LCPs against uncertainties---especially in the sense of
robust optimization---is still in its infancy.
During the last years, robust LCPs have only been studied using the
notions of strict and $\Gamma$-robustness. Unfortunately, both concepts
lead to the problem that the existence of robust solutions cannot be
guaranteed.
In this paper, we consider affinely adjustable robust LCPs. In the latter, a
part of the LCP solution is allowed to adjust via a function that is
affine in the uncertainty.
We show that this notion of robustness allows to establish strong
characterizations of solutions for the cases of uncertain matrix and
vector, separately, from which existence results can be derived.
\rev{Our main results are valid for the case of an uncertain LCP vector.
  Here, we additionally provide sufficient
  conditions on the LCP matrix for the uniqueness of a solution.
  Moreover, based on characterizations of the affinely adjustable robust
  solutions, we derive a mixed-integer programming formulation
  that allows to solve the corresponding robust counterpart.
  If, in addition, the certain LCP matrix is positive semidefinite,
  we prove polynomial-time solvability and uniqueness of robust
  solutions.}
If the LCP matrix is uncertain, characterizations of solutions are
developed for every nominal matrix, \ie, these characterizations are,
in particular, independent of the definiteness of the nominal matrix.
Robust solutions are also shown to be unique for positive definite LCP
matrix but both uniqueness and mixed-integer programming
formulations still remain open problems if the nominal LCP matrix is
not positive definite.


\end{abstract}

\keywords{Linear Complementarity Problems,
Adjustable Robustness,
Robust Optimization,
Existence,
Uniqueness%
%
}
\subjclass[2010]{90C33, 
91B50, 
91A10, 
90Cxx, 
90C34
%
%
}

\maketitle

\section{Introduction}
\label{sec:introduction}

Linear complementarity problems (LCPs) are an important tool both in
mathematical theory as well as in applied mathematics.
On the one hand, they serve as a bridge between mathematical fields
such as optimization, game theory, and matrix theory---on
the other hand, they provide one of the main modeling concepts for
market equilibrium problems in energy applications like power or gas
networks.
For an overview of these connections, we refer to the seminal
textbook \cite{Cottle_et_al:2009}.
Most likely, its strongest connection can be drawn to quadratic
programming (QP) via the fact that the Karush--Kuhn--Tucker (KKT)
conditions of many QPs can be represented as LCPs, which is also the
key aspect for the applicability of LCPs in contexts such as energy
markets; see, \eg,
\cite{Hobbs_Helman:2004,Gabriel_et_al:2012,Hobbs:1998,Metzler_et_al:2003}.

One very active sub-area of mathematical optimization in the last
decades was and is optimization under uncertainty, \ie, the study of
optimization problems in which all or a certain number of parameters
of the model are unknown or subject to perturbations.
In order to hedge against uncertainties, two major approaches
have been established:
stochastic optimization (see, \eg,
\cite{Kall_Wallace:1994,Birge_Louveaux:2011}) and robust optimization
(see, \eg, \cite{Soyster:1973,Ben-Tal_et_al:2009,Bertsimas_et_al:2011}).
While the former assumes knowledge about the distributions of the
uncertain parameters and considers, \eg, the maximization of expected
returns or the minimization of expected costs, the latter makes no
distributional assumptions but protects against the worst-case
uncertainty realization within a prescribed uncertainty set.

Although the relation between LCPs and optimization is pretty close,
comparably few research papers focus on LCPs under uncertainty.
Most of the related papers tackle the stochastic case and consider the
minimization of the expected residual gap function of the LCP; see, \eg,
\cite{Chen_Fukushima:2005,Chen_et_al:2012,Chen_et_al:2009,Lin_Fukushima:2006}
and the references therein.
In contrast to stochastic LCPs, the robust treatment of LCPs under
uncertainty is still in its infancy.
To the best of our knowledge, the first paper on robust LCPs is
\cite{Wu_et_al:2011}, in which the authors consider strict
robustifications of LCPs.
The same concept has been studied
in~\cite{Xie_Shanbhag:2014,Xie_Shanbhag:2016}.
In these contributions, the authors consider strictly robust
counterparts of uncertain LCPs for the case of different uncertainty
sets such as box or ellipsoidal uncertainties.
In particular, these papers focus on tractability of the corresponding
robust counterparts.
The results are applied to the case of
Cournot--Bertrand equilibria in power networks in~\cite{Mather_Munsing:2017}; see also
\rev{\cite{Kramer_et_al:2018,Celebi_et_al:2020} for related studies} of
  Nash--Cournot and perfect competition equilibria in comparable settings.

The concept of strict robustness in optimization has received
criticism due to the high degree of conservatism of the solutions that
it may deliver.
Consequently, several less conservative notions of robustness have been
developed during the last twenty years; see, \eg,
\cite{Sim:2004,Bertsimas_Sim:2003,Bertsimas_Sim:2004}
for $\Gamma$-robustness,
\cite{Fischetti_Monaci:2009} for light robustness,
\cite{Ben-Tal_et_al:2004,Ben-Tal_et_al:2009,Yanikoglu_et_al:2019} for
adjustable robustness, or \cite{Assmann_et_al:2018} for deciding
robustness in a fully adjustable setting with an empty first stage.
Following the idea of studying less conservative notions of
robustness, the concept of $\Gamma$-robustness has been applied to
LCPs in \cite{Krebs_Schmidt:2019} for the case of $\ell_1$-and
box-uncertainty sets and in~\cite{Krebs_et_al:2019} for the case of
ellipsoidal uncertainties.
Applications of $\Gamma$-robust LCPs in the area of power markets or
traffic equilibrium problems can be found
in~\cite{Kramer_et_al:2018,Celebi_et_al:2020,Krebs_et_al:2019}.
To the best of our knowledge, the given and rather short list of
papers on robust LCPs is complete.

Besides the study of algorithms for their solution, the most classic
topic regarding LCPs is the consideration of characterizations,
existence, and uniqueness of solutions.
These topics closely link the field to the
area of matrix classes in applied linear algebra; see again
\cite{Cottle_et_al:2009} and the many references therein.
Unfortunately, almost all the papers on robust LCPs cited above make
the observation that strong characterizations and, thus, existence of
robust solutions to LCPs cannot be ensured because the requirement
that a point is a complementarity solution for all realizations of
uncertainty is very strong.
This observation is made in~\cite{Xie_Shanbhag:2016} for strict
robustness and in~\cite{Krebs_Schmidt:2019,Krebs_et_al:2019} for
$\Gamma$-robustness.
As a remedy, the authors study the LCP's quadratic gap function and
consider the existence and uniqueness of solutions or the tractability of
problems in which the complementarity condition is not strictly
demanded but in which its violation is penalized in the LCP's gap
function.
Thus, there is one major gap in the existing literature on robust
LCPs, namely:
\begin{quote}
  \emph{Is there a robustification concept that (i)
    allows to derive strong characterizations of solutions of the
    uncertain LCP itself---instead of the LCP's gap formulation---and
    that (ii) allows to establish non-trivial robust solutions of an
    uncertain LCP?}
\end{quote}

To the best of our knowledge, only the concepts of strict and
$\Gamma$-robustness have been studied for robust LCPs.
Both do not satisfy the conditions in the question above.

In order to cure this, it is necessary to go beyond single-stage
robustness concepts---in particular, to go to two-stage robust models.
\rev{Thus, i}n this paper, we carry over the concept of adjustable
robustness to the field of LCPs under box uncertainty.
The main rationale of doing so is that the split of variables into
here-and-now as well as wait-and-see variables that can be
adjusted to the uncertainty indeed allows to characterize robust LCP
solutions and to establish non-trivial solutions.
\rev{In adjustable-robust optimization, one usually first needs to
  specify the class of functions that can be used to adjust the
  wait-and-see variables in dependence of the uncertainty.
  The easiest functions to tackle are affine functions.
  Although this may be a rather restrictive choice, it is a natural
  modelling approach that can lead to algorithmically tractable robust
  counterparts~\cite{Ouorou2011,Poss_Raack2011,Almaraj_Trafalis2020}.
  It already gives us enough flexibility to derive
  strong characterizations of robust LCP solutions as well as
  existence results.
  Thus, we first focus on affine adjustability in this paper and
  postpone more complicated uncertainty-dependent decision rules to
  our future work.}
The class of adjustable robust LCPs is introduced in
Section~\ref{sec:problem-statement} and an illustrating example is
given in Section~\ref{sec:adjustb-robust-mark}.
Afterward, we consider the cases of uncertain LCP vector and LCP
matrix separately.
\rev{Our main results are given in Section~\ref{sec:unertainty-q}
  for the case of uncertain LCP vector.
  We} derive strong characterizations of
robust solutions, from which an existence result is derived.
The used characterizations do not require any further
assumptions on the LCP matrix.
This holds both for the case of full- and lower-dimensional
uncertainty sets.
Moreover, we illustrate exemplarily the existence of non-trivial
robust LCP solutions.
Uniqueness of solutions is shown for the case of positive
(semi-)definite LCP matrix, in which we also obtain
polynomial-time solvability.
We additionally present a mixed-integer
programming formulation that can be used to compute affinely
adjustable robust LCP solutions by using standard solvers.
Characterizations of solutions can also be derived in the case of
uncertain LCP matrix; see Section~\ref{sec:uncertainty-M}.
Here, uniqueness and tractability are shown for the case of positive
definite nominal LCP matrix, whereas both remain open problems for arbitrary matrices.
The paper closes with some concluding remarks and a brief discussion
of possible topics of future work in Section~\ref{sec:conclusion}.


\section{Problem Statement}
\label{sec:problem-statement}

Given a matrix $M \in \R^{n\times n}$ and a vector $q\in\R^n$, the
linear complementarity problem LCP($q, M$) is the problem to find a
vector $z \in \R^n$ satisfying the conditions
\begin{equation}
  \label{eq:LCPold}
  z\geq 0,
  \quad Mz+q\geq 0,
  \quad z^\top(Mz+q)=0
\end{equation}
or to show that no such vector exists.
In the following, we use the standard $\perp$-notation and abbreviate
\eqref{eq:LCPold} as
\begin{equation}
  \label{eq:LCP}
  0 \leq z \perp Mz + q \geq 0.
\end{equation}
In 
real-world applications, the parameters~$M$ and~$q$ may
be uncertain.
In order to model 
this, 
we define uncertainty sets~$\Unc_M \subseteq
\R^{k_1}$ as well as $\Unc_q \subseteq
\R^{k_2}$ with suitable~$k_1$ and~$k_2$. We then consider
$M(\zeta)$ and $q(u)$ with $\zeta \in \Unc_M$ and $u \in
\Unc_q$.
The specific definition of the uncertainty sets will be
given in the corresponding sections.
Since these definitions will be qualitatively different for $M$ and
$q$ we choose to use a Greek letter to parameterize~$M$ and a Latin
letter to parameterize~$q$.

We follow the robust paradigm for dealing with such uncertain
parameters. In the strictly robust model, we want to find a
vector~$z \in \R^n$ that fulfills the conditions in~\eqref{eq:LCP} for
every possible realization of uncertainty $(\zeta, u) \in \Unc_M
\times \Unc_q$, \ie,
\begin{equation*}
  0 \leq z \perp M(\zeta)z + q(u) \geq 0
  \quad
  \fa
  \quad
  (\zeta, u) \in \Unc_M \times \Unc_q.
\end{equation*}

We call such a vector~$z$ a strictly robust solution of the uncertain LCP.
This approach is discussed in
\cite{Xie_Shanbhag:2014,Xie_Shanbhag:2016}.
The $\Gamma$-robust approach is
discussed in \cite{Krebs_Schmidt:2019,Krebs_et_al:2019}.
The main conceptual problem with strictly as well as $\Gamma$-robust
LCPs is that one usually cannot prove the existence of a solution.

The goal of this paper is to study the well established and typically less
conservative approach of (affinely) adjustable robustness in the
context of LCPs.
For adjustable robustness, a part of the solution is allowed to
adapt to a given realization of uncertainty.
The task thus is to find a vector~$r \in \R^n$,
which can be adjusted for all
uncertainties $(\zeta, u) \in \Unc_M \times \Unc_q$ by a vector
$y(\zeta, u)$ so that $z(\zeta, u) \define r + y(\zeta, u)$
satisfies
\begin{equation}
  \label{eq:ULCP-general}
  0\leq z(\zeta, u) \perp M(\zeta) z(\zeta, u) + q(u) \geq 0
  \quad
  \fa
  \quad
  (\zeta, u) \in \Unc_M \times \Unc_q.
\end{equation}
We call such a point~$z(\zeta,u)$ an adjustable robust solution of the
uncertain LCP.
In many applications, further restrictions need to be imposed on the
adjustable solution.
For instance, one usually has to distinguish between adjustable and
non-adjustable, or \enquote{here-and-now}, variables.
To this end, we introduce a parameter~$h \in \set{0, \dotsc, n}$
and require that the first $h$~entries of $y(\zeta, u)$ are zero.
This means that the first $h$~entries are non-adjustable here-and-now
decisions.

In general, the adjustable robust approach without further
assumptions on the adaptability leads to intractable problems; see,
\eg, \cite{Ben-Tal_et_al:2004}, where this is shown for the easiest
possible case of uncertain linear programs.
In this paper, we impose an assumption that is often used in adjustable robustness.
Namely, we restrict
ourselves to consider affinely adjustable robust solutions, \ie,
we restrict the solutions to be of the form
\begin{equation*}
  z(\zeta, u) = D_1\zeta + D_2u + r
  \quad \text{with} \quad
  D_1\in\R^{n\times k_1}, \,
  D_2\in\R^{n\times k_2}, \,
  r\in\R^n.
\end{equation*}
We call an affine function $z(\zeta, u) = D_1\zeta + D_2u + r$ solving
Problem~\eqref{eq:ULCP-general} an affinely adjustable robust (AAR)
solution of the uncertain LCP.
Hence, we search for affine decision rules given by~$D_1$, $D_2$, and $r$
that specify how to react to a given realization of uncertainty.
To model $h$ here-and-now variables, we w.l.o.g.\ require that the
first $h$~rows of $D_1$ and~$D_2$ are zero.

We close this section by brief\/ly introducing some notation that is
required in the remainder of this paper.
Let $A \in \R^{n\times n}$, $b \in \R^n$, and index sets~$I, J \subseteq [n]
\define \set{1, \dotsc, n}$ be given.
Then, $A_{I,J} \in \R^{|I| \times |J|}$ denotes the submatrix of~$A$
consisting of the rows indexed by~$I$ and the columns indexed by~$J$.
Moreover, $b_I$ denotes the subvector with components specified by
entries in~$I$.
If $I = J$, we also write $A_I$ instead of $A_{I,I}$. For $i,j\in[n]$
let $\delta_{ij}$ be the Kronecker delta, i.e., $\delta_{ij}=1$ if
$i=j$ and $\delta_{ij}=0$ otherwise.
Finally, the identity matrix of size~$k \times k$ is denoted
by~$\I_k$.


\section{Illustrating Example: Adjustable Robust Energy Market
  Equilibrium Modeling}
\label{sec:adjustb-robust-mark}

In this section, we consider a stylized energy market equilibrium
problem to illustrate the applicability of adjustable robustness in a
practically relevant application of market modeling.
To this end, we start with a simple market model based on the one
given in~\cite{Cottle_et_al:2009} and we also follow the notation
used there.
First, let the production sector of our energy market model be given
by the linear program
\begin{subequations}
  \label{eq:example-production}
  \begin{align}
    \min_{x \in \R^n} \quad
    & c^\top x
    \\
    \st \quad
    & Ax \geq b,
      \label{eq:example-production:technology}\\
    & Bx \geq r^*,
      \label{eq:example-production:demand} \\
    & x \geq 0,
  \end{align}
\end{subequations}
with vectors~$c \in \R^n$, $b \in \R^m$, $r^* \in \R^k$ as well as
matrices $A \in \R^{m \times n}$ and $B \in \R^{k \times n}$.
The variable vector $x$ models production levels that should be
cost-minimal but that also need to satisfy certain technological
constraints~\eqref{eq:example-production:technology} and demand satisfaction
constraints \eqref{eq:example-production:demand}.
The demand~$r^*$ itself depends on market prices~$p^*$, which is
modeled by an affine demand function, \ie,
\begin{equation*}
  r^* = D p^* + d
  \quad \text{with} \quad
  D \in \R^{k \times k}, \, d \in \R^k.
\end{equation*}
\rev{In many applications, one assumes that the matrix~$D$ is negative
  semi-definite to model that demand is non-increasing in dependence
  of the prices.
  For diagonal matrices~$D$ this then leads to monotonically decreasing and univariate demand functions, which is a very classic economic
  modeling.}
As usual in standard micro-economic settings, we need the additional
equilibrating condition $p^* = \pi^*$ with $\pi^*$ being the optimal
dual multiplier of the demand
constraint~\eqref{eq:example-production:demand}.
By using this condition as well as the (necessary and sufficient)
Karush--Kuhn--Tucker (KKT) conditions of
Problem~\eqref{eq:example-production}, the market equilibrium can be
modeled using the LCP
\begin{align*}
  0 \leq x
  & \perp c - A^\top \lambda - B^\top p^* \geq 0,
  \\
  0 \leq \lambda
  & \perp -b + Ax \geq 0,
  \\
  0 \leq \rev{p^*}
  & \perp -D p^* - d + Bx \geq 0,
\end{align*}
which is obtained by simplifying the KKT complementarity conditions
and solving for $r^*$ and $\pi^*$. The dual multiplier of the technology
constraint~\eqref{eq:example-production:technology} is denoted
by~$\lambda$.
The corresponding LCP data is given by
\begin{equation*}
  z =
  \begin{pmatrix}
    x \\ \lambda \\ p
  \end{pmatrix},
  \quad
  M =
  \begin{bmatrix}
    0 & -A^\top & -B^\top \\
    A & 0 & 0 \\
    B & 0 & -D
  \end{bmatrix},
  \quad
  q =
  \begin{pmatrix}
    c \\ -b \\ -d
  \end{pmatrix}.
\end{equation*}
If this rather general market equilibrium problem is considered as an
abstract setting for an energy market, adjustable robustness in the
context of LCPs shows up rather naturally.
Here, the electricity demand~$r^*$ depends on prices but also has a
price-insensitive part~$d$.
This vector can, for instance, be estimated 
from
historical data.
However, the demand parameter~$d$ is uncertain due to, e.g., unknown
future weather conditions, which leads to an uncertain LCP vector~$q =
q(u)$ with $u$ in some properly chosen uncertainty set~$\mathcal{U}_q$.
These uncertainties in demand can usually be tackled by adjustments
in production, \ie, not the \enquote{nominal} market equilibrium
production is used but production is adjusted in dependence of the
realization of demand uncertainty.
Since, on the other hand, certain generators such as wind or solar
power plants cannot be adjusted as easily as, \eg, coal power plants,
this additionally leads to a rather natural split between adjustable
and non-adjustable LCP variables.
\rev{Note that for $D$ being negative semidefinite, the bisymmetric
  matrix~$M$ is positive semidefinite.
  Thus, this practically relevant example belongs to the class of
  robust LCPs for which we present the strongest theoretical results
  in this paper---namely robust LCPs with uncertain vector~$q$ and
  positive semi-definite matrix~$M$.}

Similarly, uncertainty in the
coefficients of the (technological as well as demand satisfaction)
constraints leads to an uncertain LCP matrix, where again some part
of the solution corresponds to variables that can be adjusted,
the other to those that are non-adjustable.


\section{Uncertainty in $q$}
\label{sec:unertainty-q}

Throughout this section we assume that the matrix $M$ is fixed and not
affected by uncertainty.
For a given nominal vector $\bar{q}\in\R^n$ and an uncertainty
set $\Unc=\Unc_q\subseteq\R^n$, we define $q(u) \define \bar{q}+u$ for
every $u\in\Unc$.
The uncertain LCP \eqref{eq:ULCP-general} then reads
\begin{align}
  \label{eq:ULCP-q}
  0\leq z(u)~\bot~Mz(u)+q(u)\geq 0\quad\fa\quad u\in\Unc.
\end{align}
We are interested in determining AAR solutions of \eqref{eq:ULCP-q}
of the form $z(\zeta, u) =  z(u)=Du+r$ with $D\in\R^{n\times n}$ and $r\in\R^n$.
To this end, we consider a box uncertainty set
\begin{align*}
  \Unc \define \defset{u\in\mathbb{R}^n}{-\bar{u}_i\leq u_i\leq
  \bar{u}_i}
\end{align*}
that is, w.l.o.g., centered around zero.
Moreover, we split the index set $[n]$ into the set of uncertain
entries
\begin{align*}
  U \define & \defset{i\in[n]}{\bar{u}_i>0},
\end{align*}
and the set of certain entries
\begin{align*}
  S \define &\defset{i\in[n]}{\bar{u}_i=0},
\end{align*}
i.e., $[n]=U\cup S$.
For notational reasons we do not remove the columns in $D$
corresponding to $S$ but fix $D_{\bdt,S}=0$.

Recall that we require $D_{[h],\bdt}=0$ in an AAR solution, since the
first $h$ variables are non-adjustable.
For a given affine function $z(u)=Du+r$, we define the sets
\begin{align*}
  I & \define \defset{i\in [h]}{r_i\neq 0},\\
  J & \define \defset{i\in [n]\setminus[h]}{r_i\neq 0},\\
  K & \define \defset{i\in [n]}{r_i\neq 0} = I\cup J,\\
  N & \define \defset{i\in [n]}{r_i= 0} = [n]\setminus K.
\end{align*}
The assumption that the uncertainty is centered around zero
immediately leads to the following key observations.

\begin{obs}\label{obs:basic-observation-r}
  Let $z(u)=Du+r$ be an AAR solution of \eqref{eq:ULCP-q}.
  Then, $r$ is a solution of the nominal LCP($\bar{q},M$).
\end{obs}
\begin{obs}\label{obs:basic-observation-q}
  Let $z(u)=Du+r$ be an AAR solution of \eqref{eq:ULCP-q}.
  Since $z(u)\geq 0$ holds for all $u\in\Unc$, the inclusion
  \begin{align*}
    \Defset{i\in[n]}{D_{i,U}\neq 0}
    \subseteq J
  \end{align*}
  holds because, otherwise, there would exist an index~$i\notin J$ and
  an uncertainty $u'\in\Unc$ with $z_i(u') = D_{i,\bdt}u' < 0$.
\end{obs}

These observations and notations will be helpful to derive the results
in the following sections.

\subsection{Characterization and Existence of Solutions}
\label{sec:existence}

In this section, we show some general properties and
characterizations of AAR solutions.
In \Cref{Thm:Char-Uarbitrary}, we derive a system of equations that
has to be satisfied by every AAR solution.
This system of equations will be used to obtain more specific
characterizations under further assumptions on the uncertainty set.
Moreover, it admits an algorithmic approach to compute an AAR
solution, which is addressed in \Cref{sec:mixed-binary-optim}.

First, we prove a basic lemma that reformulates the constraints in
the uncertain LCP.

\begin{lemma}
  \label{lemma:uncq-basiclemma}
  \rev{Let $z(u)=Du+r$ and assume $D_{[h],\bdt}=0$. Then, the} function $z(u)$ is an
  AAR solution of \eqref{eq:ULCP-q} if and
  only if
  \begin{subequations}
    \label{eq:uncq-basiclemma}
    \begin{align}
      z_K(u) & \geq 0 \quad \fa \quad u\in\Unc,
               \label{eq:uncq-basiclemma1}\\
      (Mz(u)+q(u))_K & =0 \quad \fa \quad u\in\Unc,
                       \label{eq:uncq-basiclemma2}\\
      (Mz(u)+q(u))_N & \geq 0 \quad \fa \quad u\in\Unc.
                       \label{eq:uncq-basiclemma3}
    \end{align}
  \end{subequations}
\end{lemma}
\begin{proof}
  We show that the conditions in~\eqref{eq:uncq-basiclemma} are equivalent
  to the uncertain LCP.
  By definition of $N$ and \Cref{obs:basic-observation-q}, $z_N(u)=0$
  holds for all $u\in\Unc$.
  Thus, \eqref{eq:uncq-basiclemma1} is equivalent to $z(u)\geq 0$ for
  all $u\in\Unc$.
  If $z(u)$ satisfies \eqref{eq:uncq-basiclemma2} and
  \eqref{eq:uncq-basiclemma3}, this implies $Mz(u)+q(u)\geq 0$ for all $u\in\Unc$.
  Additionally, for all $u\in\Unc$ we have
  \begin{align*}
    z(u)^\top (Mz(u)+q(u))=z_K(u)^\top (Mz(u)+q(u))_K+z_N(u)^\top (Mz(u)+q(u))_N=0,
  \end{align*}
  where the last equality is due to \eqref{eq:uncq-basiclemma2} and
  $z_N(u)=0$.
  Thus, $z(u)$ satisfies~\eqref{eq:ULCP-q}.

  It remains to show that \eqref{eq:uncq-basiclemma2} is a necessary
  condition.
  To this end, let $z(u)=Du+r$ be an AAR solution.
  As noted before,
  $z_K(u)\geq 0$ holds for all $u\in\Unc$.
  Let us now assume that there is an index~$i\in K$
  such that there exists $u'\in\Unc$ with $z_i(u')=0$.
  This implies that $u'$ minimizes $z_i(u)=D_{i,\bdt}u+r_i$.
  Since $D_{i,\bdt}u+r_i$ is an affine function in $u$, the minimum is
  attained at the boundaries, i.e.,
  \begin{align*}
    u_j'= \begin{cases}
      \bar{u}_j, & \text{ if }D_{i,j}<0,\\
      -\bar{u}_j, & \text{ if }D_{i,j}>0,
    \end{cases}
  \end{align*}
  for all $j$ with $D_{i,j}\neq 0$.
  As $r_K>0$, we obtain $z_K(u)>0$ for all $u$ contained in the
  relative interior $\mathrm{relint}(\Unc)$.
  Furthermore, the uncertain LCP conditions imply $(Mz(u)+q(u))_K=0$ for
  all $u\in\mathrm{relint}(\Unc)$, which immediately yields
  \eqref{eq:uncq-basiclemma2} since $(Mz(u)+q(u))_K$ is an affine
  function in $u$ as well.
\end{proof}
In the following, we use Condition~\eqref{eq:uncq-basiclemma2} to
derive characterizations and properties of AAR solutions.
In \Cref{Thm:Char-Uarbitrary}, we reformulate the LCP
conditions and obtain a system of equations that needs to be satisfied by $D$ and $r$.

\begin{lemma}
  \label{Thm:Char-Uarbitrary}
  The function $z(u)=Du+r$ satisfies \eqref{eq:uncq-basiclemma2} if and
  only if $D$ and $r$ satisfy the system of equations
  \begin{subequations}
    \label{eq:uncq-char}
    \begin{align}
      M_{K\cap S,J}D_{J,U}&=0, \label{eq:uncq-char-D-S}
      \\
      M_{K\cap U,J}D_{J,K\cap U}&=-\I_{K\cap U},\label{eq:uncq-char-D-U1}\\
      M_{K\cap U,J}D_{J,N\cap U}&=0,\label{eq:uncq-char-D-U2}\\
      M_Kr_K&=-\bar{q}_K.\label{eq:uncq-char-q}
    \end{align}
  \end{subequations}
\end{lemma}
\begin{proof}
  Let $i\in K$.
  We show that $(Mz(u)+q(u))_i=0$ holds for all $u\in\Unc$ if and only
  if \eqref{eq:uncq-char} are satisfied.
  We have
  \begin{align*}
    (Mz(u)+q(u))_i&=M_{i,\bdt}z(u)+q_i(u) \\
                  &=M_{i,\bdt}Du+M_{i,\bdt}r+\bar{q}_i+u_i\\
                  &=M_{i,J}D_{J,\bdt}u+M_{i,\bdt}r+\bar{q}_i+u_i,
  \end{align*}
  where the last equality follows from $D_{i,\bdt}=0$ for all $i\notin
  J$ by \Cref{obs:basic-observation-q}.
  If $i\in K\cap S$,
  we have $u_i=0$ and thus
  \begin{equation*}
    (Mz(u)+q(u))_i=M_{i,J}D_{J,\bdt}u+M_{i,\bdt}r+\bar{q}_i=0
  \end{equation*}
  holds for all $u\in\Unc$ if and only if
  \begin{align*}
    M_{i,J}D_{J,U}=0 \quad \text{and}\quad M_{i,\bdt}r=-\bar{q}_i.
  \end{align*}
  If $i\in
  K\cap U$, $$(Mz(u)+q(u))_i=M_{i,J}D_{J,\bdt}u+M_{i,\bdt}r+\bar{q}_i+u_i=0$$
  holds for all $u\in\Unc$ if and only if
  \begin{align*}
    M_{i,J}D_{J,j} & = -\delta_{ij} \quad \fa j\in U,
    \\
    M_{i,\bdt}r & = -\bar{q}_i.\qedhere
  \end{align*}
\end{proof}

If the uncertainty set $\Unc$ is full-dimensional, i.e., $S=\emptyset$,
the system of equations~\eqref{eq:uncq-char} is rich enough to derive a
complete characterization of an AAR solution as we will show in the
following.
To this end, we first assume $S\subseteq[h]$, meaning that only the entries
of $q(u)$ corresponding to the non-adjustable variables might be certain.
Thus, the entries of $q(u)$ corresponding to adjustable variables are
all uncertain.
Under this assumption, we derive conditions that are equivalent to
\eqref{eq:uncq-char-D-S}--\eqref{eq:uncq-char-D-U2} in the following
lemma.

\begin{lemma}
  \label{Thm:Char-OnAdjFullDim}
  Let $S\subseteq [h]$.
  Then, $D$ and $r$ satisfy
  \eqref{eq:uncq-char-D-S}--\eqref{eq:uncq-char-D-U2} if and only if
  they satisfy the conditions
  \begin{subequations}
    \label{eq:Char-OnAdjFullDim}
    \begin{align}
      D_J & = -M_J^{-1},\label{eq:Char-OnAdjFullDim1}\\
      D_{J,i} & = 0 \quad \fa i\in N\cap U,\label{eq:Char-OnAdjFullDim2}\\
      I\cap U & = \emptyset,\label{eq:Char-OnAdjFullDim3}\\
      M_{I,J} & = 0.\label{eq:Char-OnAdjFullDim4}
    \end{align}
  \end{subequations}
\end{lemma}
\begin{proof}
  We first note that $S\subseteq [h]$ implies $J\subseteq U$ and thus
  $J\subseteq K\cap U$.
  Let $D$ and~$r$ satisfy \eqref{eq:uncq-char-D-S}--\eqref{eq:uncq-char-D-U2}.
  We show that they satisfy
  \eqref{eq:Char-OnAdjFullDim1}--\eqref{eq:Char-OnAdjFullDim4}.
  Since $J\subseteq K\cap U$, \eqref{eq:uncq-char-D-U1} implies
  $M_{J}D_J=-\I_J$ and thus $D_J=-M_J^{-1}$, which is
  \eqref{eq:Char-OnAdjFullDim1}.
  Furthermore, \eqref{eq:uncq-char-D-U1}
  and~\eqref{eq:uncq-char-D-U2} imply $M_JD_{J,i}=0$ for all $i\in
  (I\cup N)\cap U$.
  Since $M_J$ has full rank, it follows $D_{J,i}=0$ for
  all $i\in (I\cup N)\cap U$ and thus \eqref{eq:Char-OnAdjFullDim2} holds
  as well.
  To show \eqref{eq:Char-OnAdjFullDim3}, we assume that there exists
  an $i\in I\cap U$.
  However, $I\subseteq K$ and \eqref{eq:uncq-char-D-U1} imply $M_{i,J} D_{J,i}=-1$
  and thus $D_{J,i}\neq 0$, contradicting the previously proved statements.
  From $I\cap U=\emptyset$ and $S\subset [h]$ it follows $I=K\cap S$
  and thus \eqref{eq:uncq-char-D-S} implies $M_{I,J}D_{J,U}=0$.
  In particular, $M_{I,J}D_J=0$ holds.
  Since $\rk(D_J)=|J|$, from $M_{I,J}D_J=0$ we obtain $M_{I,J}=0$ and
  thus \eqref{eq:Char-OnAdjFullDim4}.

  Now, let $D$ and $r$ satisfy
  \eqref{eq:Char-OnAdjFullDim1}--\eqref{eq:Char-OnAdjFullDim4}.
  By direct
  insertion, it is easy to verify that~$D$ and~$r$ satisfy
  \eqref{eq:uncq-char-D-S}--\eqref{eq:uncq-char-D-U2}.
\end{proof}

We can now combine \Cref{Thm:Char-OnAdjFullDim} and
Condition~\eqref{eq:uncq-char-q} in \Cref{Thm:Char-Uarbitrary} to
obtain the desired results for the case of full-dimensional
uncertainty sets, i.e., for $S=\emptyset$.
The first one states that all non-adjustable variables
necessarily need to have a value of zero.

\begin{corollary}
  Let $S=\emptyset$ and suppose that $z(u)=Du+r$ is an AAR solution
  of~\eqref{eq:ULCP-q}.
  Then, all non-adjustable variables are zero, i.e.,
  $I=\emptyset$ and $K=J$.
\end{corollary}

Moreover, we can use the characterizations of $D$ and
$r$ from \Cref{Thm:Char-Uarbitrary} and \Cref{Thm:Char-OnAdjFullDim}
to obtain a complete characterization of AAR solutions for the case of
full-dimensional uncertainty sets.

\begin{theorem}\label{Cor:Uncq-UFullDim-Char}
  Let $S=\emptyset$.
  Then, $z(u)=Du+r$ is an AAR solution of \eqref{eq:ULCP-q} if and only
  if $D$ and $r$ are given by
  \begin{align*}
    &D_J=-(M_J)^{-1},\quad D_{J,i}=0, \, \rev{D_{i,\bdt}=0}\text{ for all } i\notin J,\\
    &r_J=-(M_J)^{-1}\bar{q}_J,
    \quad
    {r_i}=0 \text{ for all } i\notin J
  \end{align*}
  and if the following conditions are fulfilled:
  \begin{enumerate}
  \item $M_J$ is invertible,
  \item $-(M_J)^{-1}q_J(u)\geq 0$ for all $u\in\Unc$,
  \item $-M_{N,J}(M_J)^{-1}q_J(u)+q_N(u)\geq 0$ for all $u\in\Unc$.
  \end{enumerate}
\end{theorem}
The last theorem establishes a one-to-one correspondence
between an AAR solution and the set of indices of nonzero variables $J$.
Hence, to compute an AAR solution, it suffices to find a
set $J$ that fulfills the conditions (a)--(c) of the theorem.
Moreover, this characterization also allows to establish a finite and compact
existence result for AAR solutions.
\begin{corollary}\label{corollary:existence}
  Let $S=\emptyset$. For every $J\subseteq [n]\setminus[h]$, for which $M_J$ is invertible, we define
   \begin{equation*}
     A_{i,j}^J \define -|(M_J^{-1})_{i,j}\bar{u}_j|
  \end{equation*}
  for all $i,j\in J$ and
  \begin{equation*}
     C_{i,j}^J \define \rev{ -|M_{i,\bdt}(M_J^{-1})_{\bdt,j}\bar{u}_j|}
   \end{equation*}
   for all $i\in N, j\in J$.
  If there exists a subset $J\subseteq [n]\setminus[h]$ such that $M_J$ is invertible and
  \begin{align*}
    \sum_{j\in J}A_{J,j}^J - (M_J^{-1})\bar{q}_J&\geq 0,\\
    \sum_{j\in J}C_{N,j}^J \rev{- \bar{u}_N} - M_{N,J}(M_J^{-1})\bar{q}_J + \bar{q}_N &\geq 0,
  \end{align*}
  holds, then there exists an AAR solution.
\end{corollary}

The uniqueness, however, of an AAR solution is not given in general
as shown in the following example, which also illustrates the
existence of non-trivial AAR solutions.
\begin{example}
  \label{Example1}
  Consider the uncertain LCP with parameters
  \begin{align*}
    M=
    \begin{bmatrix}
      4&10\\1&2
    \end{bmatrix},
    \quad
    \bar{q}=
    \begin{pmatrix}
      -100\\-22
    \end{pmatrix},
    \quad
    \Unc=[-1,1]^2,
    \quad
    h = 0.
  \end{align*}
  There are two different AAR solutions corresponding to different
  index sets.
  For $J_1=\{1\}$, we obtain
  \begin{equation*}
    D=\begin{bmatrix}
      \rev{-\frac{1}{4}}&0\\0&0
    \end{bmatrix},
    \quad
    r=\begin{pmatrix}
      25\\0
    \end{pmatrix}
  \end{equation*}
  and for $J_2=\{2\}$, we have
  \begin{equation*}
    D=
    \begin{bmatrix}
      0&0\\0&\rev{-\frac{1}{2}}
    \end{bmatrix},
    \quad
    r=\begin{pmatrix}
      0\\11
    \end{pmatrix}.
  \end{equation*}

  Note that the matrix~$M$ is not positive semidefinite.
  We later show in \Cref{sec:positive-semidefinite-M} that being
  positive semidefinite is a sufficient condition for an
  AAR solution to be unique in the case of $S=\emptyset$.
\end{example}

\subsection{A Mixed-Integer Programming Formulation}
\label{sec:mixed-binary-optim}

In this section we make use of the reformulations given in
\Cref{Thm:Char-Uarbitrary} and state a mixed-integer feasibility
problem with binary variables that can be used to compute an AAR
solution of the uncertain LCP~\eqref{eq:ULCP-q}.
\begin{theorem}\label{Prop:Uncq-MIP}
  Let $B\in\R$ be sufficiently large and
  consider the mixed-integer feasibility problem
  \begingroup
  \allowdisplaybreaks
  \begin{subequations}\label{Eq-Uncq-MIP}
    \begin{align}
      \text{Find}~~~&x\in\{0,1\}^n,  r\in\R^n,  A,C,D\in\R^{n\times n}\\
      \st~~~
                    &Bx_i\geq r_i\geq 0,&i\in [n],\label{Eq:Uncq-MIP-r1}\\
                    &B(1-x_i)\geq M_{i,\bdt}r+\bar{q}_i\geq 0, &i\in [n],\label{Eq:Uncq-MIP-r2}\\
                    &D_{[h],\bdt}=0,~ D_{\bdt,S}=0,\label{Eq:Uncq-MIP-D1}\\
                    &\rev{B(1-x_i)\geq M_{i,\bdt}D_{\bdt, j}\geq-B(1-x_i)},&i\in S, \, j\in U\label{Eq:Uncq-MIP-D3}\\
                    &B(1-x_j)-1\geq M_{j,\bdt}D_{\bdt, j}\geq -B(1-x_j)-1,&j\in U,\label{Eq:Uncq-MIP-D4}\\
                    &B(1-x_i)\geq M_{i,\bdt}D_{\bdt, j}\geq -B(1-x_i),&i\neq j\in U,\label{Eq:Uncq-MIP-D5}\\
                    &A_{i,j}\leq -D_{i,j}\bar{u}_j, & i\in [n],\,
                                                           j\in U,\label{Eq:Uncq-MIP-A1}\\
                    &A_{i,j}\leq D_{i,j}\bar{u}_j, & i\in [n],\, j\in U,\label{Eq:Uncq-MIP-A2}\\
                    &\sum_{j\in U}A_{i,j}+r_i\geq 0,&i\in [n],\label{Eq:Uncq-MIP-A3}\\
                    &C_{i,j}\leq -(M_{i,\bdt}D_{\bdt,j}+\delta_{ij})\bar{u}_j,&i\in [n],\, j\in U,\label{Eq:Uncq-MIP-C1}\\
                    &C_{i,j}\leq (M_{i,\bdt}D_{\bdt,j}+\delta_{ij})\bar{u}_j,&i\in [n],\, j\in U,\label{Eq:Uncq-MIP-C2}\\
                    &\sum_{j\in U}C_{i,j}+M_{i,\bdt}r+\bar{q}_i\geq 0,&i\in [n].\label{Eq:Uncq-MIP-C3}
    \end{align}
  \end{subequations}
  \endgroup
  If \eqref{Eq-Uncq-MIP} is feasible, it returns an AAR solution of the form $z(u)=Du+r$ to \eqref{eq:ULCP-q}.
  If it is infeasible, then no AAR solution exists.
\end{theorem}
\begin{proof}
  It suffices to show that every solution of \eqref{Eq-Uncq-MIP} corresponds to an
  AAR solution and vice versa.
  First, let $(x,r,A,C,D)$ be a solution of \eqref{Eq-Uncq-MIP}.
  Note that $D$ fulfills the basic requirements
  $D_{[h],\bdt}=0$ and $D_{\bdt,S}=0$ by \eqref{Eq:Uncq-MIP-D1}.
  We now show that $z(u)=Du+r$ is an AAR solution.
  The inequality
  \begin{equation*}
    \sum_{j\in U}A_{i,j} \leq \min_{u\in\Unc} \Set{D_{i,\bdt}u}
  \end{equation*}
  holds for all $i\in [n]$ by \eqref{Eq:Uncq-MIP-A1} and
  \eqref{Eq:Uncq-MIP-A2}.
  It follows
  \begin{align*}
    \min_{u\in\Unc} \, \set{z_i(u)} \geq \sum_{j\in
    U}A_{i,j}+r_i \geq 0
  \end{align*}
  for all $i\in [n]$,
  where the last inequality follows from \eqref{Eq:Uncq-MIP-A3}.
  This implies $z(u)\geq 0$ for all $u\in\Unc$.
  In particular, since $r_N=0$, we also obtain $D_{N,U}=0$ and hence
  $z_N(u)=0$ for all $u\in\Unc$.
  Due to \eqref{Eq:Uncq-MIP-r1}, we have $x_i=1$ if  $i\in K$.
  Thus, \eqref{Eq:Uncq-MIP-r2} implies \eqref{eq:uncq-char-q} and
  \eqref{Eq:Uncq-MIP-D3}--\eqref{Eq:Uncq-MIP-D5} imply the
  conditions \eqref{eq:uncq-char-D-S}--\eqref{eq:uncq-char-D-U2}.
  Hence, \eqref{eq:uncq-basiclemma2} holds due to
  \Cref{Thm:Char-Uarbitrary}, i.e., $(Mz(u)+q(u))_K=0$ for all
  $u\in\Unc$.
  From $z_N(u)=0$ and $(Mz(u)+q(u))_K=0$ for all $u\in\Unc$ it
  immediately follows $z(u)^\top (Mz(u)+q(u))=0$ for all $u\in\Unc$.\\
  It remains to show that $(Mz(u)+q(u))_N\geq 0$ holds for all
  $u\in\Unc$.
  The inequalities \eqref{Eq:Uncq-MIP-C1} and \eqref{Eq:Uncq-MIP-C2}
  imply
  \begin{equation*}
    C_{i,j} \leq \min_{u\in\Unc} \Set{M_{i,\bdt}D_{\bdt,j}u_j+\delta_{ij}u_j}
  \end{equation*}
  for all $i\in [n]$, $j\in U$.
  Hence, we obtain
  \begin{equation*}
    \sum_{j\in U} C_{i,j}\leq \min_{u\in\Unc} \Set{M_{i,\bdt}Du+u_i}
  \end{equation*}
  for all $i\in N \subseteq [n]$.
  It follows
  \begin{align*}
    \min_{u\in\Unc} \Set{(Mz(u)+q(u))_i} \geq \sum_{j\in
    U}C_{i,j}+M_{i,\bdt}r+\bar{q}_i \geq 0
  \end{align*}
  for all $i\in N$,
  where the last inequality follows from \eqref{Eq:Uncq-MIP-C3}.
  Thus, $(Mz(u)+q(u))_N\geq 0$ holds for all $u\in\Unc$.

  Now, let $z(u)=Du+r$ be an AAR solution of \eqref{eq:ULCP-q}.
  Next, we construct $x$, $A$, and $C$ such that $(x,r,A,C,D)$ is a solution
  of \eqref{Eq-Uncq-MIP}.
  For all $i\in K$, we set $x_i=1$ and for all $i\in N$ we set $x_i=0$.
  Since $r$ is a nominal solution, the constraints
  \eqref{Eq:Uncq-MIP-r1} and \eqref{Eq:Uncq-MIP-r2} are satisfied for
  sufficiently large~$B$.
  Since $D$ fulfills the basic requirements $D_{h,\bdt}=0$ and
  $D_{\bdt, S}=0$, Condition~\eqref{Eq:Uncq-MIP-D1} is satisfied.
  Furthermore, $D$ is a solution of the equations
  \eqref{eq:uncq-char-D-S}--\eqref{eq:uncq-char-D-U2} in
  \Cref{Thm:Char-Uarbitrary} and, thus, $D$ satisfies
  \eqref{Eq:Uncq-MIP-D3}--\eqref{Eq:Uncq-MIP-D5} for sufficiently
  large~$B$.
  Next, we define \rev{$A_{i,j} \define -|D_{i,j}\bar{u}_j|$} for all $i,j\in[n]$.
  Then, \eqref{Eq:Uncq-MIP-A1} and \eqref{Eq:Uncq-MIP-A2} are
  satisfied, implying
  \begin{align*}
    \sum_{j\in U}A_{i,j}+r_i = \min_{u\in\Unc} \Set{z_i(u)} \geq 0
  \end{align*}
  for all $i\in [n]$.
  Hence, \eqref{Eq:Uncq-MIP-A3} is satisfied.
  Lastly, we define
  \begin{equation*}
    C_{i,j} \define -|(M_{i,\bdt}D_{\bdt,j}+\delta_{ij})\bar{u}_j|
  \end{equation*}
  for all $i,j\in[n]$.
  Then, \eqref{Eq:Uncq-MIP-C1} and \eqref{Eq:Uncq-MIP-C2} are
  satisfied, implying
  \begin{align*}
    \sum_{j\in U}C_{i,j} + M_{i,\bdt}r + \bar{q}_i
    = \min_{u\in\Unc} \Set{(Mz(u)+q(u))_i} \geq 0
  \end{align*}
  for all $i\in [n]$.
  Hence, \eqref{Eq:Uncq-MIP-C3} is satisfied.
\end{proof}

\begin{remark}
  One crucial aspect regarding the correctness of the binary
  feasibility problem in \Cref{Prop:Uncq-MIP} is that the constant~$B$ needs
  to be sufficiently large.
  For general LCPs, it can be computationally expensive to compute this
  constant; see, \eg, \cite{Pardalos:1988}.
  However, for specific instances, problem-specific structure can
  often be exploited to obtain such constants\rev{; see, \eg,
  \cite{Kleinert_Schmidt:2019a}, where similar constants are derived
  by using the specific structure of a market equilibrium
  problem that can also be modeled as a complementarity problem.}
\end{remark}

\subsection{Positive Semidefinite $M$}
\label{sec:positive-semidefinite-M}

In the remainder of this section, we assume that the matrix~$M$ is
positive semidefinite.
In this case, we attain polynomial-time solvability and
uniqueness results under further assumptions on the uncertainty
set~$\Unc$.
First, we review the following well established theorem on
linear complementarity problems.
\begin{lemma}[Theorem 3.1.7 (a), (c) in
  \cite{Cottle_et_al:2009}]\label{Lemma:Cottle-PosSemDef-Partc}
  Let $M\in\R^{n\times n}$ be positive semidefinite and let $q\in\R^n$
  be chosen arbitrarily.
  Then, the following assertions hold.
  \begin{itemize}
  \item[(a)] If $z^1$ and $z^2$ are two solutions of the LCP($q,M$), then
    \begin{equation*}
      (z^1)^\top (q+Mz^2)=(z^2)^\top (q+Mz^1)=0.
    \end{equation*}
  \item[(b)] If the LCP($q,M$) has a solution, then the set SOL($q,M$)
    of solutions is polyhedral and given by
    \begin{align*}
      \text{SOL}(q,M) = \{z \in \R^n_{\geq 0}
      & : q+Mz\geq 0,
        q^\top(z-\bar{z}) = 0,
      \\
      & \quad (M+M^\top)(z-\bar{z}) = 0\},
    \end{align*}
    where $\bar{z}$ is an arbitrary solution.
  \end{itemize}
\end{lemma}

For what follows, we define
\begin{align*}
  P \define \Defset{j\in [n]}{\exists z\in
  \text{SOL}(\bar{q},M):z_j>0},
  \quad
  L \define [n] \setminus P.
\end{align*}
\rev{For the following results, we need to know the index set~$P$
  explicitly.
  Note that $\text{SOL}(\bar{q},M)$ can be explicitly stated via
  Part~(b) of the previous lemma since the special solution~$\bar{z}$
  can be computed by solving a single convex quadratic program.
  The set~$P$ can then be obtained by solving $n$ linear programs in which
  $z_j$, $j \in [n]$, is maximized over the polyhedral feasible
  set~$\text{SOL}(\bar{q},M)$ and by checking afterward, whether the
  solution is strictly positive.
  Thus, $P$ can be computed in polynomial time.}

We now use \Cref{Lemma:Cottle-PosSemDef-Partc} to strengthen
\Cref{Thm:Char-Uarbitrary}.

\begin{lemma}
  \label{lemma:Uncq-Char-Mpossemdef}
  Let $M$ be positive semidefinite.
  If $z(u)=Du+r$ is an AAR solution of \eqref{eq:ULCP-q}, the system
  of equations
  \begin{subequations}
    \begin{align}
      M_{P\cap S,P}D_{P,U}&=0\label{eq:Uncq-Char-Mpossemdef-PS},\\
      M_{P\cap U,P}D_{P,P\cap U}&=-\I_{P\cap U},\label{eq:Uncq-Char-Mpossemdef-PUPU}\\
      M_{P\cap U,P}D_{P,L\cap U}&=0.\label{eq:Uncq-Char-Mpossemdef-PULU}
    \end{align}
  \end{subequations}
  is satisfied.
\end{lemma}
\begin{proof}
  From \Cref{obs:basic-observation-r} we know that $r$ is a nominal
  solution.
  Thus, due to \Cref{Lemma:Cottle-PosSemDef-Partc} (a),
  $M_{P,\bdt}r+\bar{q}_P=0$ holds.
  Since $z(u)=Du+r$ is an AAR solution, we know
  \begin{equation*}
    Mz(u)+q(u)=MDu+Mr+\bar{q}+u\geq 0
  \end{equation*}
  for all $u\in\Unc$.
  In particular, we have
  \begin{align*}
    (MDu+Mr+\bar{q}+u)_P =
    M_{P,\bdt}Du + M_{P,\bdt}r + \bar{q}_P + u_P
    = M_{P,\bdt}Du + u_P
    \geq 0
  \end{align*}
  for all $u\in\Unc$.
  Since we set $D_{\bdt,S}=0$, we have $D_{P,\bdt}u=D_{P,U}u_U$ and
  from \Cref{obs:basic-observation-q} it follows
  $M_{P,\bdt} D = M_P D_{P,\bdt}$.
  Hence, the inequality
  \begin{align*}
    M_{P,P}D_{P,U}u_U+u_P\geq 0
  \end{align*}
  holds for all $u\in\Unc$.

  For $i\in P\cap S$, we have $u_i=0$ and, thus, $M_{i,P}D_{P,U}u_U\geq 0$
  holds for all $u\in \Unc$.
  This implies $M_{i,P}D_{P,U}=0$ as otherwise there would
  exist an element~$u'\in\Unc$ from the uncertainty set defined by
  $u'_U=-\lambda(M_{i,P}D_{P,U})^\top $ for some $\lambda>0$ and
  $u'_S=0$ so that $M_{i,P}D_{P,U}u'_U=-\lambda
  ||M_{i,P}D_{P,U}||_2<0$.
  Thus, \eqref{eq:Uncq-Char-Mpossemdef-PS} holds.

  Next, for $i\in P\cap U$ we have $M_{i,P}D_{P,U}u_U+u_i\geq 0$
  for all $u\in\Unc$.
  For the same reasons as in the previous case, this implies
  $M_{i,P}D_{P,U}u_U=-u_i$, as otherwise we could again construct an
  uncertainty $u'$ in the box uncertainty set $\Unc$ so that
  $M_{i,P}D_{P,U}u'_U+u'_i < 0$.
  We obtain \eqref{eq:Uncq-Char-Mpossemdef-PUPU} and
  \eqref{eq:Uncq-Char-Mpossemdef-PULU}.
\end{proof}

We now combine \Cref{Lemma:Cottle-PosSemDef-Partc} and
\ref{lemma:Uncq-Char-Mpossemdef} to obtain a linear feasibility
problem that can be used to solve the uncertain LCP with positive
semidefinite~$M$.
Thus, in this case, there is no need to solve the mixed-integer
feasibility problem from \Cref{Prop:Uncq-MIP}.
\begin{theorem}
  \label{Thm:Uncq-Mpsd-FeasProb}
  Let $M$ be positive semidefinite and suppose further that $\bar{z}$
  is a solution of the nominal LCP($\bar{q},M$).
  Consider the linear feasibility problem
  \begingroup
  \allowdisplaybreaks
  \begin{subequations}
    \label{Eq:Uncq-Mpsd-FeasProb}
    \begin{align}
      \text{Find} \quad &r\in\R^n, \,  A,C,D\in\R^{n\times n}\nonumber\\
      \st \quad &\rev{r\geq 0,} ~\bar{q}+Mr\geq 0,\label{Eq:Uncq-Mpsd-FeasProb-r1}\\
                    &\bar{q}^\top(r-\bar{z})=0,\label{Eq:Uncq-Mpsd-FeasProb-r2}\\
                    &(M+M^\top)(r-\bar{z})=0,\label{Eq:Uncq-Mpsd-FeasProb-r3}\\
                    &D_{L,\bdt}=0,~ D_{[h],\bdt}=0,~ D_{\bdt,S}=0,\label{Eq:Uncq-Mpsd-FeasProb-D1}\\
                    &\begin{bmatrix}M_{P\cap S,P}\\M_{P\cap U,P}\end{bmatrix}
      \begin{bmatrix}D_{P,P\cap U}&D_{P,L\cap U}\end{bmatrix}
                                =\begin{bmatrix}0&0\\-\I_{P\cap U}&0\end{bmatrix},\label{Eq:Uncq-Mpsd-FeasProb-D2}\\
                    &A_{i,j}\leq -D_{i,j}\bar{u}_j, \quad i\in P, \, j\in U,\label{Eq:Uncq-Mpsd-FeasProb-A1}\\
                    &A_{i,j}\leq D_{i,j}\bar{u}_j, \quad i\in P, \, j\in U,\label{Eq:Uncq-Mpsd-FeasProb-A2}\\
                    &\sum_{j\in U}A_{i,j}+r_i\geq 0,\quad i\in P,\label{Eq:Uncq-Mpsd-FeasProb-A3}\\
                    &C_{i,j}\leq
                      -(M_{i,\bdt}D_{\bdt,j}+\delta_{ij})\bar{u}_j,\quad
                      i\in L, \, j\in U,\label{Eq:Uncq-Mpsd-FeasProb-C1}\\
                    &C_{i,j}\leq
                      (M_{i,\bdt}D_{\bdt,j}+\delta_{ij})\bar{u}_j,\quad
                      i\in L, \, j\in U,\label{Eq:Uncq-Mpsd-FeasProb-C2}\\
                    &\sum_{j\in U}C_{i,j}+M_{i,\bdt}r+\bar{q}_i\geq 0,
                      \quad i\in L.\label{Eq:Uncq-Mpsd-FeasProb-C3}
    \end{align}
  \end{subequations}
  \endgroup
  Every feasible point of \eqref{Eq:Uncq-Mpsd-FeasProb} corresponds to
  an AAR solution of the form $z(u)=Du+r$.
  If \eqref{Eq:Uncq-Mpsd-FeasProb} is infeasible, then no AAR solution
  exists.
\end{theorem}
As parts of the proof of \Cref{Thm:Uncq-Mpsd-FeasProb} are similar to that
of \Cref{Prop:Uncq-MIP}, we keep the following proof rather short.
\begin{proof}
  Let $(r,A,C,D)$ be a solution of \eqref{Eq:Uncq-Mpsd-FeasProb}.
  We show, that $z(u)=Du+r$ is an AAR solution.
  First, we note that $D$ satisfies the basic requirements
  $D_{[h],\bdt}=0$ and $D_{\bdt,S}=0$ by~\eqref{Eq:Uncq-Mpsd-FeasProb-D1}.
  Since $r$ satisfies
  \eqref{Eq:Uncq-Mpsd-FeasProb-r1}--\eqref{Eq:Uncq-Mpsd-FeasProb-r3},
  it is a solution of the nominal LCP($\bar{q},M)$ by
  \Cref{Lemma:Cottle-PosSemDef-Partc} (b).
  Therefore, we obtain $r_L=0$ by the definition of $P$ and $L$,
  $D_{L,\bdt}=0$ by \eqref{Eq:Uncq-Mpsd-FeasProb-D1} and thus
  $z_L(u)=0$ holds for all $u\in\Unc$.
  Furthermore, we know $(Mr+\bar{q})_P=0$ due to
  \Cref{Lemma:Cottle-PosSemDef-Partc} (a).
  From \eqref{Eq:Uncq-Mpsd-FeasProb-D2} it follows $(MDu)_P=-u_P$ and
  thus
  \begin{equation*}
    (Mz(u)+q(u))_P=(MDu)_P+u_P+(Mr+\bar{q})_P=0
  \end{equation*}
  holds for all $u\in\Unc$.
  From $z_L(u)=0$ and $(Mz(u)+q(u))_P=0$ for all $u\in \Unc$ it
  follows $z(u)^\top (Mz(u)+q(u))=0$ for all $u\in\Unc$.

  It remains to show that $z_P(u)\geq 0$ and $(Mz(u)+q(u))_L\geq 0$
  holds for all $u\in\Unc$.
  The constraints
  \eqref{Eq:Uncq-Mpsd-FeasProb-A1}--\eqref{Eq:Uncq-Mpsd-FeasProb-A3}
  imply $z_P(u)\geq 0$ for all $u\in\Unc$ and the constraints
  \eqref{Eq:Uncq-Mpsd-FeasProb-C1}--\eqref{Eq:Uncq-Mpsd-FeasProb-C3}
  imply  $(Mz(u)+q(u))_L\geq 0$ for all $u\in\Unc$ for the same
  reasons as in the proof of \Cref{Prop:Uncq-MIP}.

  Now, let $z(u)=Du+r$ be an AAR solution.
  We construct $A$ and $C$ such that $(r,A,C,D)$ is a solution of
  \eqref{Eq:Uncq-Mpsd-FeasProb}.
  We know that $r$ is a nominal solution and, thus,
  \eqref{Eq:Uncq-Mpsd-FeasProb-r1}--\eqref{Eq:Uncq-Mpsd-FeasProb-r3}
  are satisfied by \Cref{Lemma:Cottle-PosSemDef-Partc}(b).
  By definition of $P$ and $L$, we have $r_L=0$ and thus $D_{L,U}=0$
  due to \Cref{obs:basic-observation-q}.
  The requirements $D_{h,\bdt}=0$ and $D_{\bdt,S}=0$ hold by definition.
  Hence, \eqref{Eq:Uncq-Mpsd-FeasProb-D1} is satisfied.
  The constraint \eqref{Eq:Uncq-Mpsd-FeasProb-D2} holds due to
  \Cref{lemma:Uncq-Char-Mpossemdef}.
  For all $i,j\in[n]$, we now define
  \begin{align*}
    A_{i,j} \define &-|D_{i,j}\bar{u}_j|,
    \\
    C_{i,j} \define &-|(M_{i,\bdt}D_{\bdt,j}+\delta_{ij})\bar{u}_j|.
  \end{align*}
  Then, \eqref{Eq:Uncq-Mpsd-FeasProb-A1}--\eqref{Eq:Uncq-Mpsd-FeasProb-C3}
  are satisfied for the same reasons as in the proof of \Cref{Prop:Uncq-MIP}.
\end{proof}

If the matrix~$M$ is positive semidefinite, the
nominal LCP can be solved by solving a convex quadratic program.
Therefore, a solution $\bar{z}$ for the nominal LCP, which we need as a
precondition in \Cref{Thm:Uncq-Mpsd-FeasProb}, can be computed in
polynomial time.
Since the linear feasibility problem
\eqref{Eq:Uncq-Mpsd-FeasProb} can be solved in polynomial time as
well, we obtain the following complexity result.

\begin{corollary}
  \label{Cor:PolyTime}
  Let $M$ be positive semidefinite.
  Then one can find an AAR solution of \eqref{eq:ULCP-q} or correctly
  state that there is no AAR solution in polynomial time.
\end{corollary}

We now use \Cref{lemma:Uncq-Char-Mpossemdef} to obtain uniqueness
results under additional assumptions on the uncertainty set.
As in the general case in \Cref{sec:existence}, we first consider
the case $S\subseteq [h]$.

\begin{lemma}
  \label{Thm:Uncq-OnAdjFullDim-Uniqueness}
  Let $M$ be positive semidefinite and $S\subseteq [h]$.
  If $z(u)=Du+r$ is an AAR solution of \eqref{eq:ULCP-q}, the matrix
  $D$ is uniquely determined by $D_{P\cap U}=-(M_{P\cap U})^{-1}$ and
  $D_{i,j} = 0$ for all $i,j \notin P\cap U$.
\end{lemma}
\begin{proof}
  From \Cref{lemma:Uncq-Char-Mpossemdef} we know
  \begin{align*}
    \begin{bmatrix}
      -\I_{P\cap U}&0
    \end{bmatrix}=M_{P\cap U,P}\begin{bmatrix}
      D_{P,P\cap U}&D_{P,L\cap U}
    \end{bmatrix}.
  \end{align*}
  Since $D_{[h],\bdt}=0$ and $S\subseteq [h]$ holds, we have
  $D_{P\cap S,\bdt}=0$, which implies
  \begin{align*}
    \begin{bmatrix}
      -\I_{P\cap U}&0
    \end{bmatrix}
    =M_{P\cap U,P}\begin{bmatrix}
      D_{P,P\cap U}&D_{P,L\cap U}
    \end{bmatrix}=M_{P\cap U}\begin{bmatrix}
      D_{P\cap U}&D_{P\cap U,L\cap U}
    \end{bmatrix}.
  \end{align*}
  Thus, the equation $M_{P\cap U}D_{P\cap U}=-\I_{P\cap U}$ implies
  $D_{P\cap U}=-(M_{P\cap U})^{-1}$.
  Furthermore, since $M_{P\cap U}$ is invertible and
  $M_{P\cap U}D_{P\cap U,L\cap U}=0$ holds, it follows $D_{P\cap U,L\cap U}=0$.
  As $D_{i,\bdt}=0$ for all $i\notin P$ due to
  \Cref{obs:basic-observation-q}, this finishes the proof.
\end{proof}

The previous lemma asserts the uniqueness of the matrix~$D$.
If we now assume that all entries of $q(u)$ are uncertain, i.e.,
$S = \emptyset$, \Cref{Thm:Uncq-OnAdjFullDim-Uniqueness} leads to 
uniqueness of the entire AAR solution.

\begin{theorem}
  \label{Thm:Uncq-UFullDim-UniqueSol}
  Let $M$ be positive semidefinite and $S = \emptyset$.
  \begin{enumerate}
  \item If there are multiple solutions to the nominal
    LCP($\bar{q},M$), there is no AAR solution.
  \item If there exists an AAR solution, it is unique.
  \end{enumerate}
\end{theorem}
\begin{proof}
  We first note that $P\cap U=P$ holds since $S = \emptyset$.
  Any solution $r$ to the nominal LCP($\bar{q},M$) satisfies
  $M_Pr_P=-\bar{q}_P$ due to \Cref{Lemma:Cottle-PosSemDef-Partc}(a)
  and the definition of~$P$.
  If there are multiple solutions, $M_P$ cannot be invertible and, thus,
  there cannot exist an AAR solution according to
  \Cref{Thm:Uncq-OnAdjFullDim-Uniqueness}.
  Hence, if there is an AAR solution $z(u)=Du+r$, $r$ is unique due to
  the previous argument and $D$ is unique due to
  \Cref{Thm:Uncq-OnAdjFullDim-Uniqueness}.
\end{proof}

We close this section with some remarks on the connection between
our results and the classical LCP theory as well as on the limits of
affine adjustability.
If the matrix~$M$ is positive semidefinite, the nominal LCP can be
solved by solving a convex QP, which can be done in polynomial time.
This is also the underlying reason for our complexity result
\Cref{Cor:PolyTime}.
As for nominal LCPs, uniqueness of solutions cannot be
guaranteed in the case of an arbitrary matrix~$M$.
\rev{Under the assumption that $M$ is a $P$ matrix, \ie, all principal
  minors of~$M$ are positive, the uniqueness of the solution of the
  nominal LCP is guaranteed for every~$q$; see, \eg,
  \cite[Chapter~3]{Cottle_et_al:2009}.}
This statement directly carries over to uncertain LCPs
with general uncertainty sets. If the solution $z(u)$ for every
realization of the uncertainties $u\in\Unc$ is unique, an AAR solution
is unique as well.
However, \Cref{Thm:Uncq-UFullDim-UniqueSol} states that, in the case
of full-dimensional uncertainty sets, we only need
positive semidefiniteness of the matrix $M$ to guarantee the
uniqueness of an AAR solution, which is a less strong condition than
$M$ being a $P$ matrix.

Note that we illustrated the existence of non-trivial solutions, see
\Cref{Example1}, and stated conditions for the existence of a
solution in \Cref{corollary:existence}.
However, let us also note that there exist uncertain LCPs that have
an adjustable but not an affinely adjustable robust
solution as the following example shows.
\begin{example}
  Consider the uncertain LCP given by
  \begin{align*}
    M=\begin{bmatrix}
      1&\frac{1}{2}\\\frac{1}{2}&1
    \end{bmatrix}\rev{\succ\,} 0,
    \quad \bar{q} =
    \begin{pmatrix}
      -5\\-3
    \end{pmatrix}, \quad \Unc=[-1,1]^2, \quad h=0.
  \end{align*}
  \rev{Since all principal minors of $M$ are positive, $M$ is a $P$ matrix.
  Hence, for any realization $u'\in \Unc$, there exists a solution of the nominal
  LCP($q(u'),M$). Therefore, a fully adjustable solution would map every realization
  to its respective unique solution.}
  However, the uncertain LCP does not have an AAR solution, which can be verified by
  applying \Cref{Cor:Uncq-UFullDim-Char}.
\end{example}

Solving the uncertain LCP with other decision rules than affine ones
is left for future research.


\section{Uncertainty in $M$}
\label{sec:uncertainty-M}

In this section, we assume that the vector~$q$ is certain and consider
uncertainty only in the matrix~$M$. In particular, we are given matrices
$M^0, M^1, \dots, M^k\in\R^{n\times n}$ as well as $\Unc = \Unc_M =
[-1,1]^k$ and define
\begin{equation*}
  M(\zeta) \define M^0 + \sum_{i=1}^k \zeta_i M^i.
\end{equation*}
The uncertain LCP~\eqref{eq:ULCP-general} then reads
\begin{align}
  \label{eq:ULCP-M}
  0\leq z(\zeta) \perp M(\zeta) z(\zeta) + q \geq 0
  \quad
  \fa
  \quad
  \zeta\in\Unc.
\end{align}
For this problem, we are interested in computing an AAR solution of
the form $z(\zeta)=D\zeta+r$ with $D\in\R^{n\times k}$ and
$r\in\R^n$. As before, we assume that the first $h$~rows of $D$ are
zero for some fixed $h$ to distinguish between adjustable and
non-adjustable variables.
However, the results presented in this section are independent of the
specific choice of $h$.
\begin{remark}
  We can interpret $M^0$ as the nominal matrix that is perturbed by
  the matrices $M^1, \dotsc ,M^k$.
  This definition of a matrix uncertainty set is considered
  in~\cite{Xie_Shanbhag:2014} for the first time in the context of
  LCPs and is also used in~\cite{Krebs_Schmidt:2019}.
\end{remark}

For an AAR solution $z(\zeta)=D\zeta + r$ we define the sets
\begin{align*}
  J \define \Defset{j\in[n]}{r_j> 0},
  \quad
  N \define [n] \setminus J.
\end{align*}
As in \Cref{obs:basic-observation-q} for the case of uncertain $q$,
we have $\defset{j\in[n]}{D_{j,\bdt}\neq 0} \subseteq J$ and, thus,
$D_{N,\bdt} = 0$.
Analogously to the proof of \Cref{lemma:uncq-basiclemma},
we have $z_J(\zeta)>0$ for all $\zeta\in\mathrm{int}(\Unc)$.

We now prove necessary conditions that every AAR solution satisfies.
\begin{theorem}
  \label{prop:uncM-basicprop}
  Let $z(\zeta)=D\zeta+r$ be an AAR solution for \eqref{eq:ULCP-M}.
  Then,
  \begin{subequations}
    \begin{align}
      M^0_Jr_J+q_J & = 0,\label{eq:uncM-basicprop1}\\
      M^i_Jr_J+M^0_JD_{J,i} & =0 \quad \fa \quad  i\in [k],\label{eq:uncM-basicprop2}\\
      M^i_JD_{J,i} & = 0 \quad \fa \quad i\in [k],\label{eq:uncM-basicprop3}\\
      M^i_JD_{J,j}+M^j_JD_{J,i} & = 0 \quad \fa \quad i,j\in [k], \, i\neq
      j \label{eq:uncM-basicprop4}
    \end{align}
  \end{subequations}
  holds.
\end{theorem}
\begin{proof}
  Since $0 \in \mathrm{\Unc}$, the vector $z(0) = r$ is a solution of the
  nominal LCP($q, M^0$) and thus \eqref{eq:uncM-basicprop1} holds.
  For $i\in[k]$, we define
  \begin{equation*}
    \Unc_i \define \Defset{\zeta\in\Unc}{\zeta_i\in (-1,1), \,
      \zeta_j=0 \text{ for all } j\neq i}
    \subseteq \mathrm{int}(\Unc).
  \end{equation*}
  We have $z_J(\zeta) > 0$ for all $\zeta\in\Unc_i$ and thus
  $(M(\zeta)z(\zeta)+q)_J=0$ holds for all $\zeta \in \Unc_i$.
  We obtain
  \begin{align*}
    0 & = (M^0_J+\zeta_iM^i_J)(r_J+\zeta_iD_{J,i})+q_J
    \\
      & = M^0_Jr_J+\zeta_i\left(M^i_Jr_J +
        M^0_JD_{J,i}\right)+\zeta_i^2M^i_JD_{J,i}+q_J
  \end{align*}
  for all $\zeta_i\in(-1,1)$.
  Hence, the conditions~\eqref{eq:uncM-basicprop2}
  and~\eqref{eq:uncM-basicprop3} follow.

  Now, for $i,j\in[k]$ with $i\neq j$, we define
  \begin{equation*}
    \Unc_{i,j} \define
    \Defset{\zeta\in\Unc}{\zeta_i, \zeta_j \in (-1,1), \, \zeta_p=0
      \text{ for all } p \notin \set{i, j}}
    \subseteq \mathrm{int}(\Unc).
  \end{equation*}
  As before, $z_J(\zeta) > 0$ holds for all $\zeta \in \Unc_{i,j}$
  and thus
  \begin{align*}
    0&=(M(\zeta)z(\zeta)+q)_J\\
     &=M(\zeta)_J(D\zeta+r)_J+q_J\\
     &=M^0_J(\zeta_iD_{J, i}+\zeta_jD_{J,j}+r_J) +
       \zeta_iM^i_J(\zeta_iD_{J, i}+\zeta_jD_{J,j}+r_J)\\
     & \quad +\zeta_jM^j_J(\zeta_iD_{J, i}+\zeta_jD_{J,j}+r_J)+q_J\\
     &=(M^0_Jr_J+q_J)
       +\zeta_i(M^0_JD_{J,i}+M^i_Jr_J)
       +\zeta_j(M^0_JD_{J,j}+M^j_Jr_J)\\
     & \quad +\zeta_i^2M^i_JD_{J,i}
       +\zeta_j^2M^j_JD_{J,j}+\zeta_i\zeta_j(M^i_JD_{J,j}+M^j_JD_{J,i})=(*)
  \end{align*}
  for all $\zeta\in\Unc_{i,j}$.
  The first term is zero due to \eqref{eq:uncM-basicprop1}.
  Applying \eqref{eq:uncM-basicprop2} and
  \eqref{eq:uncM-basicprop3}, all other terms except for the last one
  are zero as well.
  It follows
  \begin{align*}
    0 = (*) = \zeta_i \zeta_j(M^i_JD_{J,j}+M^j_JD_{J,i})
  \end{align*}
  for all $\zeta\in\Unc_{i,j}$ and thus \eqref{eq:uncM-basicprop4}
  holds.
\end{proof}

Since the systems of equations of the last theorem might allow for
multiple solutions, they are not sufficient to fully characterize an
AAR solution.
However, under the additional assumption that $M^0_J$ is invertible,
it is possible to derive a complete characterization.
For example, this assumption is satisfied if $M^0$ is positive
definite as in this case every submatrix $M^0_I$, $I \subseteq [n]$,
is invertible.

We first introduce some notation and subsequently present the
complete characterization in \Cref{prop:UncM-M0pd-Char}.
To this end, let $M^0_I$ be invertible for a subset $I \subseteq [n]$.
Then, we define
\begin{equation*}
  \tM^{I,i} \define (M^0_I)^{-1} M^i_I (M^0_I)^{-1}.
\end{equation*}

\begin{corollary}
  \label{prop:UncM-M0pd-Char}
  Let $z(\zeta)=D\zeta+r$ be an AAR solution for \eqref{eq:ULCP-M}.
  If $M_J^0$ is invertible, then $D$ and $r$ are given by
  \begin{align*}
    D_{J,i} = \tM^{J,i} q_J, \, i\in[k],
    \quad
    r_J = -(M^0_J)^{-1}q_J,
    \quad
    D_{N,\bdt}=0,
    \quad
    r_N=0.
  \end{align*}
\end{corollary}
\begin{proof}
  Since $M_J^0$ is invertible, \eqref{eq:uncM-basicprop1} is
  equivalent to $r_J=-(M^0_J)^{-1}q_J$.
  By using this equation for $q_J$,
  \eqref{eq:uncM-basicprop2} can be equivalently reformulated as
  $M^0_JD_{J,i}=M^i_J(M^0_J)^{-1}q_J$ for all $i\in [k]$.
  Thus, for all $i\in [k]$ we obtain
  \begin{equation*}
    D_{J,i}=(M^0_J)^{-1}M^i_J(M^0_J)^{-1}q_J=\tM^{J,i} q_J.
    \qedhere
  \end{equation*}
\end{proof}

In the next example, we illustrate that indeed solutions
characterized by this corollary exist.
\begin{example}
  Let
  \begin{equation*}
    M(\zeta)=\begin{bmatrix}
      4&1\\0&4
    \end{bmatrix}+\zeta
              \begin{bmatrix}
                0&1\\0&0
              \end{bmatrix}, \quad
                        q=\begin{pmatrix}
                          -8\\-16
                        \end{pmatrix}, \quad h=0.
  \end{equation*}
  As $M^0$ is invertible, we consider the set $J = [n]$.
  It follows
  \begin{align*}
    (M^0)^{-1}=\frac{1}{16}\begin{bmatrix}
      4&-1\\0&4
    \end{bmatrix}\quad\text{and}\quad
               \tM^{J,1}=(M^0)^{-1}M^1(M^0)^{-1}=\frac{1}{16}\begin{bmatrix}
                 0&1\\0&0
               \end{bmatrix}.
  \end{align*}
  Using \Cref{prop:UncM-M0pd-Char}, we obtain
  \begin{align*}
    r=-(M^0)^{-1}q=\begin{pmatrix}
      1\\4
    \end{pmatrix},
    \quad
    D=\tM^{J,1} q=\begin{pmatrix}
      -1\\0
    \end{pmatrix}.
  \end{align*}
  It is easy to verify that $z(\zeta)=(1-\zeta, 4)^\top$ is an AAR
  solution.
\end{example}

For what follows, let $z(\zeta)=D\zeta +r$ be an AAR solution and
suppose that $M_J^0$ is invertible.
The conditions~\eqref{eq:uncM-basicprop3} and~\eqref{eq:uncM-basicprop4} can be
reformulated similarly as in the proof of \Cref{prop:UncM-M0pd-Char}
by using the characterizations of $r$ and $D$.
We obtain that \eqref{eq:uncM-basicprop3} is equivalent to
\begin{equation*}
  M^i_J\tM^{J,i}q_J=0 \quad \fa \quad i\in [k].
\end{equation*}
Expression~\eqref{eq:uncM-basicprop4} is equivalent to
\begin{equation*}
  (M^i_J\tM^{J,j} + M^j_J\tM^{J,i})q_J = 0
  \quad \fa \quad
  i,j\in [k], \, i\neq j.
\end{equation*}
We combine these conditions and obtain
\begin{equation}
  \label{eq:UncM-M0pd-condq}
  q_J\in\bigcap_{i,j\in[k]} \ker\left(M^i_J\tM^{J,j}+M^j_J\tM^{J,i}\right).
\end{equation}

In the following, we derive a reformulation of the uncertain LCP
conditions in~\eqref{eq:ULCP-M} such that they only depend on the LCP
parameters~$M$ and~$q$.
To this end, we use \Cref{prop:UncM-M0pd-Char}.
The equation
\begin{equation}
  \label{eq:UncM-M0pd-Dzeta}
  D_{J,\bdt}\zeta = \sum_{i=1}^k\zeta_i\tM^{J,i} q_J
\end{equation}
holds for all $\zeta\in\Unc$.
Thus, the requirement that $z_J(\zeta)\geq 0$ for all $\zeta\in\Unc$
is equivalent to
\begin{equation}
  \label{eq:UncM-M0pd-condMq1}
  \left( \sum_{i\in[k]} \zeta_i \tM^{J,i} - (M^0_J)^{-1}\right) q_J
  \geq 0
  \quad \fa \quad
  \zeta\in\Unc.
\end{equation}
Furthermore, $(M(\zeta)z(\zeta)+q)_N\geq 0$ for all $\zeta\in\Unc$ is
equivalent to
\begin{equation}
  \label{eq:UncM-M0pd-condMq2}
  M_{N,J}(\zeta)\left(\sum_{i\in[k]}\zeta_i\tM^{J,i}-(M^0_J)^{-1}\right)q_J
  + q_N \geq 0
  \quad\fa\quad
  \zeta \in \Unc.
\end{equation}
The following theorem summarizes that these conditions lead to a
full characterization.
\begin{theorem}
  \label{thm:final-char-unc-M}
  Let $D$ and $r$ be characterized as in \Cref{prop:UncM-M0pd-Char} for
  $J \subseteq [n]$ such that $M_0^J$ is invertible.
  \rev{Furthermore, suppose that $D_{[h],\bdt}=0$ holds.}
  Then, $z(\zeta)=D\zeta+r$ is an AAR solution for \eqref{eq:ULCP-M} if and
  only if $M(\zeta)$ and $q$ fulfill the
  conditions~\eqref{eq:UncM-M0pd-condq}, \eqref{eq:UncM-M0pd-condMq1},
  and \eqref{eq:UncM-M0pd-condMq2}.
\end{theorem}
\begin{proof}
  It only remains to show that $(M(\zeta)z(\zeta)+q)_J=0$ for all
  $\zeta\in\Unc$ is implied by~\eqref{eq:UncM-M0pd-condq}.
  For all $\zeta\in\Unc$ we have
  \begin{align*}
    (M(\zeta)z(\zeta)+q)_J&=\rev{M_J(\zeta)z_J(\zeta)}+q_J\\
                          &=M_J(\zeta)(D_{J,\bdt}\zeta+\rev{r_J})+q_J\\
                          &=M^0_JD_{J,\bdt}\zeta+\sum_{i=1}^k
                            \zeta_iM^i_JD_{J,\bdt}\zeta +
                            M^0_Jr_J + \sum_{i=1}^k\zeta_iM^i_Jr_J+q_J\\
                          &=M^0_JD_{J,\bdt}\zeta+\sum_{i=1}^k\zeta_iM^i_JD_{J,\bdt}\zeta
                            + \sum_{i=1}^k\zeta_iM^i_Jr_J=(*),
  \end{align*}
  where we used $M^0_Jr_J = -q_J$.
  We apply \eqref{eq:UncM-M0pd-Dzeta} and obtain
  \begin{align*}
    (*)&=M^0_J\sum_{i=1}^k \zeta_i \tM^{J,i}
         q_J+\sum_{i=1}^k\zeta_iM^i_J\sum_{j=1}^k\zeta_j\tM^{J,j}
         q_J+\sum_{i=1}^k\zeta_iM^i_Jr_J
    \\
       &=\sum_{i=1}^k \zeta_i M^i_J(M^0_J)^{-1}
         q_J+\sum_{i,j\in[k]}\zeta_i\zeta_j
         M^i_J\tM^{J,j}q_J+\sum_{i=1}^k\zeta_iM^i_Jr_J=(**).
  \end{align*}
  By \eqref{eq:UncM-M0pd-condq} we know $\sum_{i,j\in[k]}\zeta_i\zeta_j M^i_J\tM^{J,j}q_J=0$. Thus,
  \begin{align*}
    (**)&=\sum_{i=1}^k \zeta_i M^i_J(M^0_J)^{-1} q_J+\sum_{i=1}^k\zeta_iM^i_Jr_J\\
        &=\sum_{i=1}^k \zeta_i M^i_J(M^0_J)^{-1}
          q_J-\sum_{i=1}^k\zeta_iM^i_J(M^0_J)^{-1}q_J = 0.
          \qedhere
  \end{align*}
\end{proof}

We conclude this section with some final remarks on the derived
results and the uniqueness of solutions. \Cref{prop:UncM-M0pd-Char}
shows that we can fully characterize an AAR solution if the nominal
matrix $M_J^0$ is invertible. In general, the difficulty lies in
finding the set $J$ of nonzero entries in the solution. Therefore,
there might exist different AAR solutions even if $M_I^0$ is
invertible for every $I\subseteq [n]$. However, if $M^0$ is positive
definite, $r$ is unique and therefore the set $J$ is unique, yielding
the uniqueness of an AAR solution if it exists at all.

Note that we do not state a general existence result here for the case
of uncertain~$M$ as we did in \Cref{corollary:existence} for
uncertain~$q$.
We think that an analogous result can be obtained, in principle, by
using \Cref{thm:final-char-unc-M} and by checking all vertices of the
box-uncertainty set for $\zeta$ in \eqref{eq:UncM-M0pd-condMq1} and
\eqref{eq:UncM-M0pd-condMq2}.
Although finite, the number of conditions in such an existence result most
likely would be exponential in the dimension of the uncertainty set.
We think that the same also holds for the size of a corresponding
mixed-integer programming formulation, which is why we omit to state
it here.

Finally, let us also comment on the case in which both the LCP
vector~$q$ as well as the LCP matrix~$M$ are uncertain.
The easier setting then is the one in which both uncertainties are
independent.
However, already this case is rather challenging for affinely
adjustable robust LCPs.
Consider, for instance, Condition~\eqref{eq:UncM-M0pd-condq}, which is
also part of the final characterization in
\Cref{thm:final-char-unc-M}.
A simultaneous consideration of $q$ and $M$ would require that the
null-space condition in~\eqref{eq:UncM-M0pd-condq} is satisfied for
$q_J(u)$ for all $u \in \mathcal{U}_q$.
Our hypothesis is that this extended condition alone would already be
rather hard to satisfy in practically meaningful LCP settings, which
is why we postpone the consideration of uncertainty in $q$ and $M$ to
future research.


\section{Conclusion}
\label{sec:conclusion}

In this paper, we studied affinely adjustable robust linear
complementarity problems with box-uncertainties either in the LCP
matrix~$M$ or in the LCP vector~$q$.
We addressed the topics of characterization, existence, and
uniqueness of solutions completely for the case of uncertain~$q$.
Moreover, we developed a mixed-integer linear model that
allows to compute affinely adjustable robust LCP solutions with
standard solvers.
For the case of uncertain~$M$, characterizations are established as
well and uniqueness of solutions is shown under the assumption that
the nominal LCP matrix is positive definite.

While the standard single-stage modeling assumptions of
strict as well as of $\Gamma$-robustness both fail to
enable the study of robust solutions directly (instead, the LCP's gap
function formulation is usually considered), imposing the assumption
of affine adjustability in the second stage is sufficient.
Thus, adjustable robustness is the first established concept of robust
optimization that has been carried over to LCPs,
which allows for studying the robust LCP solutions directly instead of
considering the gap function formulation as a replacement.
However, several problems remain open.
For instance, a compact existence result and a compact mixed-integer
programming formulation for the case of uncertain LCP matrix is
missing.
Moreover, the consideration of other uncertainty sets like ellipsoids
or the consideration of non-affine decision rules is part of our
future research.


\section*{Acknowledgments}
\label{sec:acknowledgements}

This research has been performed as part of the Energie Campus
Nürnberg (EnCN) and is supported by funding of the Bavarian State
Government.
The authors thank the Deutsche Forschungsgemeinschaft for their
support within project A05, B06, and B08 in the
Sonderforschungsbereich/Transregio 154 \enquote{Mathematical
  Modelling, Simulation and Optimization using the Example of Gas
  Networks}.


\printbibliography

\end{document}